\documentclass[a4paper, 12pt]{article}

\usepackage[left=3cm,right=3cm,
    top=3cm,bottom=3cm,bindingoffset=0cm]{geometry}
\usepackage{graphicx}
\usepackage{amssymb,amsfonts,amsthm, amscd}
\usepackage{amsmath}

\usepackage{xcolor}
\usepackage{hyperref}

\begin{document}
\title{Automorphism tower problem and semigroup of endomorphisms for free Burnside groups}
\author{V.S.Atabekyan}

\date{}


\newtheorem{thm}{Theorem}
\newtheorem{cor}[thm]{Corollary}
\newtheorem{lem}{Lemma}
\newtheorem{prop}{Proposition}
\newtheorem{defn}{Definition}
\newtheorem{rem}{\bf{Remark}}
\newtheorem{alg}{\bf{Algorithm}}

\maketitle

\begin{abstract}
We have proved that  the group of all inner
automorphisms of the
free Burnside group $B(m,n)$ is the unique normal subgroup in
$Aut(B(m,n))$ among all its subgroups, which are isomorphic to free
Burnside group $B(s,n)$ of some rank $s$ for all odd 
$n\ge1003$ and $m>1$. It follows that the group of automorphisms $Aut(B(m,n))$ of the
free Burnside group $B(m,n)$ is complete for odd $n\ge1003$, that is it has a trivial center and any automorphism of $Aut(B(m,n))$ is inner. Thus, for groups
$B(m,n)$ is solved the automorphism tower problem and is showed that
it is as short as the automorphism tower of the absolutely free
groups. Moreover,  proved that every
automorphism of $End(B(m,n))$ is a conjugation by an element of $Aut(B(m,n))$.

\end{abstract}

\paragraph{Introduction.}
If the center of a group $G$ is trivial, then it is embedded into
the group of its automorphisms $Aut(G)$. Such embedding is given by
mapping each element of group to the inner automorphism generated by
this element. The inner automorphism generated by an element $g\in
G$ is denoted by $i_g$ and is defined by the formula
$x^{i_g}=gxg^{-1}$ for all $x\in G$ (the image of the element $x$
under the map $\alpha$ is denoted by $x^\alpha$). The easily
verifiable relation for composition of automorphisms $\alpha\circ
i_g\circ\alpha^{-1}= i_{g^\alpha}$ shows that the group of all inner
automorphisms $Inn(G)$ is a normal subgroup in $Aut(G)$. Moreover,
the relation $\alpha\circ i_g\circ\alpha^{-1}=i_{g^\alpha}$ implies
that in a group $G$ with trivial center the centralizer of the
subgroup $Inn(G)$ is also trivial in $Aut(G)$. In particular, the
group $Aut(G)$ is also a centerless group. This allows to consider
the automorphism tower
\begin{equation}\label{tow}
G=G_0\lhd G_1 \lhd\cdot\cdot\cdot\lhd G_k\lhd\cdot\cdot\cdot,
\end{equation}
where $G_k=Aut(G_{k-1})$ and $G_k$ is identified with $Inn(G_k)$
under the embedding $G_k\hookrightarrow Aut(G_k)$, $g\mapsto i_g$
($k=1,2,...$).

According to classical Wielandt's theorem (see \cite[Theorem
13.5.2]{Rob}), the automorphism tower \eqref{tow} of any finite
centerless group  terminates after a finite number of steps (that is
there exists a number $k$ such that $G_k=G_{k-1}$). For infinite
groups the analogous statement is false (for example, automorphism
tower of infinite dihedral group does not terminate in finitely many
steps). In early seventies G.~Baumslag proposed to study the
automorphism tower for absolutely free and for some relatively free
groups. In particular, he formulated a hypothesis that the tower of
absolutely free group of finite rank should be ``very short''.

\medskip
In 1975 J.~Dyer and E.~Formanek in \cite{DF} confirmed the
Baumsalag's hypothesis proving that if $F$ is a free group of finite
rank $>1$, then its group of automorphisms $Aut(F)$ is complete.
Recall that a group is called \textit{complete}, if it is centerless
and each of its automorphisms is inner. V.~Tolstikh in \cite{T}
proved the completeness of $Aut(F)$ for free groups $F$ of infinite
rank. It is clear that if the group of automorphisms of centerless
group $G_0$ is complete, then its automorphism tower terminates
after the first step, that is, $G_0\lhd G_1=G_2=...\,.$

The later new proofs and various generalizations of Dayer-Formanek
theorem have been obtained by E.~Formanek \cite{F}, D.G.~Khramtsov
\cite{Kh}, M.R.~Bridson and K.~Vogtmann \cite{BV}.

Further, in \cite{26} and \cite{24} it was established that the
group of automorphisms of each non-abelian free solvable group of
finite rank is complete. It was showed that the group of
automorphisms of free nilpotent group of class 2 and rank $r\ge 2$
is complete provided that $r \not= 3$. In the case n = 3 the height
of the automorphism tower \eqref{tow} is 2.

Note that in all above-mentioned results on automorphism tower of
relatively free groups only torsion free groups were considered.

\paragraph{\bf Preliminary and the main results.}
We study the automorphism tower of free Burnside groups $B(m,n)$,
i.e. the relatively free groups of rank $m>1$ of the variety of all
groups which satisfy the identity $x^n=1$. The group $B(m,n)$ is the
quotient group of absolutely free group $F_m$ on $m$ generators by
normal subgroup $F_{m}^n$ generated by all $n$-th powers. Obviously,
any periodic group of exponent $n$ with $m$ generators is a quotient
group of $B(m,n)$. By the theorem of S.I.~Adian, the group $B(m,n)$
of rank $m>1$ is infinite for any odd $n\ge665$. This Theorem and a
series of fundamental properties of $B(m, n)$, was proved in the
monograph \cite{A}. A comprehensive survey of results about the free
Burnside groups and related topics is given in the paper \cite{UMN}.

Our main result states that the automorphism tower of non-cyclic
free Burnside group $B(m,n)$ is terminated on the first step for any
odd $n\ge 1003$. Hence, the automorphism tower problem for groups
$B(m,n)$ is solved. We is show that it is as short as the
automorphism tower of the absolutely free groups. In particular, the
group $Aut(B(m,n))$ is complete. The inequality $n\ge 1003$ for the
exponent $n$ is closely related to the result in the last chapter of
monograph \cite{A}, concerning the construction of an infinite
independent system of group identities in two variables (the
solution of the finite basis problem). Based on the developed
technique in this chapter in \cite {AL} the authors have constructed
an infinite simple group of period $ n $ with cyclic subgroups for
each $ n \ge 1003 $, which plays a key role in our proof of the main
result. The use of simple quotient groups to obtain information
about the automorphisms of $B(m,n)$ first occurs in the paper
\cite{O03} of A.Yu.Olshanskii. We are pleased to stress the
influence of \cite{O03} on our research.

\medskip The well known Gelder-Bear's theorem asserts that every
complete group is a direct factor in any group in which it is
contained as a normal subgroup (see \cite[Theorem 13.5.7]{Rob}).
According to Adian's theorem (see \cite[Theorem 3.4]{A}) for any odd
$n\geq665$ the center of (non-cyclic) free Burnside group is
trivial. However, the groups $B(m,n)$ are not complete, because, for
example, the automorphism $\phi$ of $B(m,n)$, defined on the free
generators by the formula $\forall i (\phi(a_i)=a_i^2)$, is an outer
automorphism.

Nevertheless, the free Burnside groups possess a property analogous
to the above-mentioned  characteristic property of complete groups.
It turns out that each group $B(m, n)$ is a direct factor in every
periodic group of exponent $n$, in which it is contained as a normal
subgroup. This statement was proved for large enough odd $n$
($n>10^{80}$) by E.~Cherepanov in \cite{Cher} and for all odd $n\ge
1003$ by the author in \cite{Ata11}.

\medskip

Our main result is the following

\begin{thm}\label{2}For any odd $n\ge1003$ and $m>1$, the group of
all inner automorphisms $Inn(B(m,n))$ is the unique normal subgroup
of the group $Aut(B(m,n))$ among all subgroups, which are isomorphic
to a free Burnside group $B(s,n)$ of some rank $s$.
\end{thm}

The proof of this result is essentially based on the papers
\cite{Af09}--\cite{A13} of the author.

Theorem \ref{2} immediately implies the following

\begin{cor}
The groups of automorphisms $Aut(B(m,n))$ and $Aut(B(k,n))$ are
isomorphic if and only if $m=k$ (for any odd $n\ge1003$).
\end{cor}

By Burnside criterion, if the group of all inner automorphisms
$Inn(G)$ of a centerless group $G$ is a characteristic subgroup in
$Aut(G)$, then $Aut(G)$ is complete (see \cite[Theorem
13.5.8]{Rob}). Since the image of the subgroup $Inn(B(m,n))$ under
every automorphism of the group $Aut(B(m,n))$ is a normal subgroup,
Theorem \ref{2} implies

\begin{cor}\label{1}
The group of automorphisms $Aut(B(m,n))$ of the free Burnside group
$B(m,n)$ is complete for any odd $n\ge1003$ and $m>1$.
\end{cor}

\medskip
It should be emphasized that the group $Aut (B (m, n))$ is saturated
with a lot of subgroups, which are isomorphic to some free Burnside
group. It is known that each non-cyclic subgroup of $B(m, n)$, and
hence the group $Inn(B(m,n))$, contains a subgroup isomorphic to the
free Burnside group $B(\infty,n)$ of infinite rank (see
\cite[Theorem 1]{A09}). Furthermore, $Aut(B(m,n))$ contains  free periodic
subgroups having trivial intersection with $Inn(B(m,n))$ for $m>2$.
For instance, consider a subgroup of $Aut(B(m,n))$ generated by
automorphisms $l_{j}$, $j=2,... ,m$, defined on generators $a_i$,
$i=1,... , m,$ by formulae $ l_{j}(a_1)=a_1a_j$ and $
l_{j}(a_k)=a_k$ for $k=2,... , m$. It is easy to check the equality
$$W(l_{2},...,l_{m})(a_1)= a_1 \cdot W(a_2,...,a_m)$$ for any word
$W(a_2,...,a_m)$. Then, the automorphism $W(l_{2},...,l_{m})$ is the
identity automorphism if and only if the equality $W(a_2,...,a_m)=1$
holds in $B(m,n)$. Hence, the automorphisms $l_{1j}$, $j=2,... ,m$
generate a subgroup isomorphic to free Burnside group $B(m-1,n)$ of
rank $m-1$. According to Theorem \ref{2}, none of the
above-mentioned subgroups is a normal subgroup of $Aut(B(m,n))$.

\medskip
A detailed proof of the above-mentioned results see \cite{ijac}

\medskip
 Using the
completeness of $Aut(B(m,n))$, we have also obtained a description of the group of automorphisms $Aut(End(B(m,n)))$ for endomorphism's semigroup of $B(m,n)$. Our result is the following
\begin{thm}
If  $S: End(B(m,n)) \to End(B(m,n))$ is an automorphism of the semigroup $End(B(m,n))$, then there is an $\alpha\in Aut(B(m,n))$ such that
$S(\beta) = \alpha\circ\beta\circ\alpha^{-1}$ for all $\beta\in
End(B(m,n))$, where $n \geq 1003$ is arbitrary odd period and $m>1$.
\end{thm}

E.Formanek proved the analogus statement for absolutely free groups of finite rank (see \cite{FP}). Note that the question about the description of $Aut(End(F))$ for relatively free groups has been posed by B.Plotkin (see \cite{MP}).

\end{document}